\documentclass[11pt,reqno]{amsart}

\usepackage[T1]{fontenc}
\usepackage{lmodern}
\usepackage[margin=1in]{geometry}
\usepackage{amsmath,amssymb}
\usepackage{xcolor}
\usepackage{microtype}
\usepackage{indentfirst,comment}

\usepackage[
  unicode=true,
  colorlinks=true,
  linkcolor={blue!75!black},
  citecolor={blue!80!black},
  urlcolor={blue!65!black},
  bookmarksnumbered=true
]{hyperref}
\usepackage{bookmark}
\hypersetup{
  pdftitle={Positivity and negativity for h-bases for n},
  pdfauthor={Wouter van Doorn, Pietro Monticone, and Melvyn B. Nathanson},
  pdfsubject={Positivity and negativity for h-bases for n},
  pdfkeywords={h-bases for n, sumsets, postage stamp problem,
    additive number theory, combinatorial number theory, AI proofs},
  pdflang={en-GB}
}

\numberwithin{equation}{section}
\theoremstyle{plain}
\newtheorem{theorem}{Theorem}[section]
\newcommand{\bt}{\begin{theorem}}
\newcommand{\et}{\end{theorem}} 

\newtheorem{corollary}[theorem]{Corollary}
\newcommand{\bc}{\begin{corollary}}
\newcommand{\ec}{\end{corollary}} 

\newtheorem{lemma}[theorem]{Lemma}
\newcommand{\bl}{\begin{lemma}}
\newcommand{\el}{\end{lemma}}

\newtheorem{problem}[theorem]{Problem}
\newcommand{\bprob}{\begin{problem}}
\newcommand{\eprob}{\end{problem}}

\AddToHook{env/theorem/before}{\begin{samepage}}
\AddToHook{env/theorem/after}{\end{samepage}}
\AddToHook{env/lemma/before}{\begin{samepage}}
\AddToHook{env/lemma/after}{\end{samepage}}
\AddToHook{env/corollary/before}{\begin{samepage}}
\AddToHook{env/corollary/after}{\end{samepage}}
\AddToHook{env/problem/before}{\begin{samepage}}
\AddToHook{env/problem/after}{\end{samepage}}

\newcommand{\beq}{\begin{equation}}
\newcommand{\eeq}{\end{equation}}
\newcommand{\N}{\ensuremath{ \mathbf N }}
\newcommand{\Z}{\ensuremath{\mathbf Z}}

\newcommand{\benum}{\begin{enumerate}}
\newcommand{\eenum}{\end{enumerate}}

\DeclareMathOperator{\qqand}{\qquad\text{and}\qquad}

\title{Positivity and negativity for $h$-bases for $n$} 

\author[van Doorn]{Wouter van Doorn}
\address{Groningen, The Netherlands} 
\email{wonterman1@hotmail.com} 

\author[Monticone]{Pietro Monticone}
\address{Harmonic, London, UK}
\email{pietro.monticone@harmonic.fun}

\author[Nathanson]{Melvyn B. Nathanson}
\address{Lehman College (CUNY), Bronx, NY 10468}
\email{melvyn.nathanson@lehman.cuny.edu}

\date{\today}
 
\begin{document}


\begin{abstract}
A set $A$ of integers is an \emph{$h$-basis for $n$} if every integer in
$[0,n] = \{0,1,2,\ldots,n\}$ belongs to the sumset $hA$, that is, can be
represented as the sum of exactly $h$ not necessarily distinct elements of $A$.
The classical postage stamp problem asks for the largest integer $n_h(k)$
such that $[0,n_h(k)] \subseteq hA$ for some set $A$ of nonnegative integers
of size $k+1$. We consider the effect of allowing negative elements.
Let $n^{\flat}_h(k)$ be the largest integer such that
$[0,n^{\flat}_h(k)] \subseteq hA$ for some set $A$ of integers of size $k+1$
containing at least one negative integer. For $k=2$, we have
$n^{\flat}_h(2) < n_h(2)$ if $h$ is odd or $h=4$ 
and $n^{\flat}_h(2) = n_h(2)$ if $h$ is even and $h \neq 4$.
Proofs of these comparison results were autonomously discovered and formally
verified in Lean by
\textnormal{\href{https://aristotle.harmonic.fun/}{Aristotle}}, a formal
reasoning agent developed by
\textnormal{\href{https://www.harmonic.fun/}{Harmonic}}.
In this paper we give new proofs of these results.
\end{abstract}

\maketitle

\section{Extremal functions for bases of order \texorpdfstring{$h$}{h} for \texorpdfstring{$n$}{n}}

Let $\N_0 = \{0,1,2,3,\ldots\}$ be  the set of nonnegative integers 
and \Z\ the set of all integers.  
Define  the interval of integers $[u,v] = \{n \in \Z: u \leq n \leq v\}$.  
Let 
\[
\binom{\N_0}{k+1} = \{A \subseteq \N_0: |A| = k+1\} 
\] 
and 
\[
\binom{\Z}{k+1} = \{A \subseteq \Z: |A| = k+1\}.
\]
Then  
\[
\binom{\Z}{k+1} \setminus \binom{\N_0}{k+1}  = \{A \subseteq \Z: |A| = k+1 
\text{ and }  a <0 \text{ for some } a \in A \}. 
\]

Let $A$ be a set of integers.  
For every positive integer $h$, the \emph{$h$-fold sumset of $A$}, denoted $hA$, 
is the set of all sums of $h$ not necessarily distinct elements of $A$.   
The set $A$ is a \emph{basis of order $h$ for $n$} 
(or an $h$-basis for $n$) if $[0,n] \subseteq hA$.  
Let $n_h(A)$ be the largest integer $n$ such that $A$ is a 
basis of order $h$ for $n$, that is, 
\[
n_h(A) = \max\{ n: [0,n]\subseteq hA\}.
\]
If $0 \notin hA$, then $n_h(A)$ is undefined. 

For example, if 
\[
A = \{-1,0,2\}
\]
then 
\[
hA = [-h,2h-2] \cup \{2h\}
\]
and 
\[
n_h(A) = 2h-2.
\]

If $A$ is a set of nonnegative integers, then $n_h(A) \geq 0$ if and only if $0 \in A$ 
and $n_h(A) \geq 1$ if and only if $0,1\in A$  if and only if $n_h(A) \geq h$. 
The classical \emph{postage stamp problem} for sets of nonnegative integers  
(introduced by Hans Rohrbach~\cite{rohr37a,rohr38} in 1937) 
is to compute  the extremal function 
\[
n_h(k)  =  \max \left\{ n_h(A): A \in \binom{\N_0}{k+1}  \right\}.  
\]
An \emph{extremal basis for $n_h(k)$} is a set $A$ of $k+1$ nonnegative integers 
such that $n_h(A) = n_h(k)$.  If $A$ is an extremal basis for $n_h(2)$, 
then $A = \{0,1,c\}$ for some integer $c > 1$.  

Alfred St\"{o}hr~\cite{stoh55}  solved the postage stamp problem for $k=2$ 
in 1955. A proof is included in~\cite{nath2026-z}.  
In Theorem~\ref{Aristotle:theorem:new} we give a new proof 
of St\"ohr's formula for $n_h(2)$. 

Let $[u]$ denote the integer part of the real number $u$, that is,
the unique integer $n$ such that $n \leq u <n+1$. 

\bt[St\"{o}hr]                \label{Aristotle:theorem:Stohr}
For all positive integers $h$, 
\[
n_h(2) = \left[ \frac{h^2+6h+1}{4}\right].
\]
For odd $h$, there is one extremal basis: 
\[
  \left\{0,1 , \frac{h+3}{2}  \right\}. 
\] 
For even $h$, there are two   extremal bases: 
\[
   \left\{0,1,  \frac{h}{2} + 1  \right\} \qqand   
\left\{0,1, \frac{h}{2} + 2 \right\}.   
\]
\et

Nathanson~\cite{nath2026-z} asked how the extremal problem changes 
if we consider sets that include 
at least one negative integer, that is, how does the extremal function 
\[
n^{\flat}_h(k) = \max\left\{n_h(A): A \in\binom{\Z}{k+1} \setminus \binom{\N_0}{k+1}\right\}
\]
compare with  the extremal function $n_h(k)$? 
An \emph{extremal basis for $n^{\flat}_h(k)$} is a set $A$ of $k+1$ integers 
such that $n_h(A) = n^{\flat}_h(k)$ and there exists $a \in A$ with $a < 0$. 
To compute $n^{\flat}_h(k)$ and to construct the extremal bases for $n^{\flat}_h(k)$ 
is much harder than the analogous problems for $n_h(k)$, even for $k=2$.  

The AI system Aristotle~\cite{Aristotle} answered a ``negativity''  question
for $k=2$, that is, for sets $A = \{a, b, c\}$ 
with $a < b < c$ and $a < 0$. 
Aristotle showed that  $n^{\flat}_h(2) < n_h(2)$ for all odd $h$ and for $h = 4$, 
and that $n^{\flat}_h(2) = n_h(2)$ for all even $h \neq 4$.  
In this paper we give (revised and improved) proofs of these results. 
In a follow-up paper we aim to compute the exact value of $n^{\flat}_h(2)$ for all $h$.

\section{Three lemmas}

\bl         \label{Aristotle:lemma:3A}
Let $A = \{a < b < c\}$ be a set of integers and let $x = b-a$ and $y = c-a$.  
Then $n \in hA$ if and only if there exist nonnegative integers $s_n$ and $t_n$ such that 
\beq            \label{Aristotle:3Aineq}
s_n + t_n \leq h 
\eeq
and 
\beq            \label{Aristotle:3Asum}
n = ha + s_nx+t_ny.
\eeq
\el

\begin{proof}
We have $n \in hA$ if and only if there exist nonnegative integers 
$r_n$, $s_n$, and $t_n$ such that  
\[
r_n + s_n + t_n  = h 
\]
and 
\[
n = r_n a + s_nb+t_nc.
\]
Then $s_n + t_n \leq h $ and 
\begin{align*}
n & = r_n a + s_nb+t_nc\\
&  = (r_n + s_n + t_n)a + s_n(b-a) +t_n(c-a) \\
& = ha + s_nx +t_ny.
\end{align*}
Thus, $n \in hA$ implies~\eqref{Aristotle:3Aineq} and~\eqref{Aristotle:3Asum}. 

Conversely, if conditions~\eqref{Aristotle:3Aineq} and~\eqref{Aristotle:3Asum} are
satisfied, then 
\[
r_n = h-s_n-t_n \in \N_0
\]
 and 
\begin{align*}
n & = ha + s_nx +t_ny  \\ 
& =  (r_n + s_n + t_n)a + s_n(b-a) + t_n(c-a) \\
& = r_n a + s_nb+t_nc \in hA. 
\end{align*}
Thus,~\eqref{Aristotle:3Aineq} and~\eqref{Aristotle:3Asum} imply $n \in hA$.
This completes the proof. 
\end{proof}

\bl         \label{Aristotle:lemma:xya}
Let $A = \{a < b < c\}$ be a set of integers and let $x = b-a$ and $y = c-a$.   
If $[0,h+1] \subseteq hA$, then 
\beq                           \label{Aristotle:xya-y}
1 \leq x < y \leq h+1
\eeq
and 
\beq                           \label{Aristotle:xya-gcd}
\gcd(x,y) = 1
\eeq
and 
\beq                           \label{Aristotle:xya-a}
a > -h.
\eeq
\el

\begin{proof}  
We have $[0,h+1] \subseteq hA$.  
By Lemma~\ref{Aristotle:lemma:3A}, for all $n \in [0,h+1]$ there exist 
$s_n,t_n \in [0,h]$ such that $s_n+t_n \leq h$ and 
\[
n = ha+s_nx+t_ny.
\] 
By the pigeonhole principle, there exist $n_1,n_2$ with $0 \leq n_1 < n_2 \leq h+1$ 
and $s_{n_1} = s_{n_2}$.  Then $t_{n_1} < t_{n_2}$ and 
\[
h+1 \geq n_2-n_1 = \left( t_{n_2} - t_{n_1}\right) y \geq y.
\] 
This proves~\eqref{Aristotle:xya-y}.

We have 
\[
0 = ha + s_0 x + t_0 y \in hA
\] 
and 
\[
1 = ha + s_1 x + t_1 y \in hA.
\] 
Subtracting, we obtain  
\[
1 = (s_1 - s_0)x + (t_1 - t_0)y
\]
and so
\[
\gcd(x,y)=1.
\]
This proves~\eqref{Aristotle:xya-gcd}.

Because  $h+1\in hA$ and $y \leq h+1$, we have  
\[
h+1= ha +s_{h+1}x+t_{h+1}y \leq ha + (s_{h+1}+t_{h+1})y \leq ha + h(h+1) 
\] 
and so 
\[
a \geq -h-1+  \frac{h+1}{h} > -h,
\]
proving~\eqref{Aristotle:xya-a}. 
\end{proof}

Let $\varphi$ be the Euler phi-function.

\bl                     \label{Aristotle:lemma:algorithm}
Let $h \geq 3$.  Then $n^{\flat}_h(2)$ 
is the maximum of $n_h(A)$ for all sets $A = \{a,a+x,a+y\}$ 
such that 
\[
1 \leq x < y \leq h+1, \qquad \gcd(x,y) = 1 \qquad \text{and} \qquad -h < a < 0.
\]
The number of sets $A$ satisfying these conditions is
\[
(h-1)\sum_{y=2}^{h+1}\varphi(y) \leq (h-1)\binom{h+1}{2}.
\]
\el

\begin{proof}
Let $h \geq 3$ and $A = \{-1,0,2\}$.  Then 
\[
n^{\flat}_h(2) \geq n_h(A) \geq 2h-2 \geq h+1.   
\] 
Thus, to compute $n^{\flat}_h(2)$, it suffices to compute $n_h(A)$
only for sets $A$ such that $[0,h+1] \subseteq hA$.  
The result now follows directly from Lemma~\ref{Aristotle:lemma:xya}.  
\end{proof}

\section{A general upper bound} 

\bt             \label{Aristotle:theorem:new}
For all positive integers $h$, 
\[
n^{\flat}_h(2) \leq n_h(2) = \left[ \frac{h^2 +6h + 1}{4} \right].
\]
\et

\begin{proof} 
For even $h$, consider the sets  
\[
A_1 = \left\{0,1,\frac{h}{2}+1\right\} \qqand A_2 = \left\{0,1,\frac{h}{2}+2\right\}.
\]
The sumset $hA_1$ contains the set 
\begin{align*} 
\left[ 0, \frac{h^2+6h}{4} \right] 
& = 
\left\{ \left(\frac{h}{2}+1\right)q+r: q \in \left[0,\frac{h}{2} \right]  
\text{ and } r \in  \left[0,\frac{h}{2} \right]   \right\} \\
&  \qquad  \cup \left\{ \left(\frac{h}{2}+1\right)^2+r: 
 r \in  \left[0,\frac{h}{2} -1\right]   \right\}
\end{align*}
and 
\[
\frac{h^2+6h}{4} +1 \notin hA_1.
\]
Therefore, 
\[
n_h(A_1) = \frac{h^2 +6h}{4}.
\]

The sumset $hA_2$ contains the set 
\begin{align*} 
\left[ 0, \frac{h^2+6h}{4} \right] 
& = 
\left\{ \left(\frac{h}{2}+2\right)q+r: q \in \left[0,\frac{h}{2} - 1\right]  
\text{ and } r \in  \left[0,\frac{h}{2} +1\right]   \right\} \\
& \qquad \cup \left\{ \left(\frac{h}{2}+2\right) \frac{h}{2}+r: 
 r \in  \left[0,\frac{h}{2}  \right]   \right\}
\end{align*}
and 
\[
\frac{h^2+6h}{4} +1 \notin hA_2.
\]
Therefore, 
\[
n_h(A_2) = \frac{h^2 +6h}{4}.   
\]

For odd $h$, consider the set   
\[
A_3 = \left\{0,1,\frac{h+3}{2} \right\}.
\]
The sumset $hA_3$ contains the set 
\begin{align*} 
\left[ 0, \frac{h^2+6h+1}{4} \right] 
& = 
\left\{ \left(\frac{h+3}{2} \right)q+r: q \in \left[0,\frac{h-1}{2} \right]  
\text{ and } r \in  \left[0,\frac{h+1}{2} \right]   \right\} \\
& \qquad \cup \left\{ \left(\frac{h+3}{2} \right) \left( \frac{h+1}{2} \right) +r: 
 r \in  \left[0,\frac{h-1}{2}  \right]   \right\}
\end{align*}
and 
\[
\frac{h^2+6h+1}{4} +1 \notin hA_3.
\]
Therefore, 
\[
n_h(A_3) = \frac{h^2 +6h+1}{4}.
\]
It follows that 
\beq        \label{Aristotle:nhA-ineq}
n_h(2) \geq \left[ \frac{h^2 +6h + 1}{4} \right] 
\eeq 
for all positive integers $h$. 
We shall prove that, for all sets $A \in \binom{\Z}{3}$,   
\beq        \label{Aristotle:nhA-ineq-2}
n_h(A) \leq \left[ \frac{h^2+6h+1}{4} \right].
\eeq

Let $A = \{a,b,c\}$ with $a < b < c$ and let $x = b-a$ and $y = c-a$.  
Inequality~\eqref{Aristotle:nhA-ineq-2} holds if $n_h(A) \leq h$.  
Thus, we can assume that 
\[
N = n_h(A) \geq h+1.
\]
By Lemma~\ref{Aristotle:lemma:xya}, 
\[
y \leq h+1
\qqand
\gcd(x,y) = 1.
\]
Consider the congruence class $ha-x \pmod{y}$. There exist unique integers 
\[
n_3 \in [0,y-1] \qqand n_4 \in [N-y+1,N]
\]
such that 
\[
n_3 \equiv n_4 \equiv ha-x \pmod{y}.
\]
Then $n_3,n_4 \in hA$ and 
\beq                  \label{Aristotle:n3n4ineq}
n_4 - n_3  \geq (N-y+1) - (y-1) = N - 2y + 2. 
\eeq
\begin{samepage}
The congruence
\[
n_3 = ha+s_{n_3}x+ t_{n_3}y \equiv ha+s_{n_3}x  \equiv ha-x \pmod{y}
\]
implies  
\[
(s_{n_3} +1)x  \equiv 0 \pmod{y}.
\]
\end{samepage}
From $\gcd(x,y)=1$, we obtain  
\[
s_{n_3} \equiv -1 \pmod{y}. 
\]
Then $s_{n_3} \geq 0$ implies  
\[
s_{n_3} \geq y-1
\]
and 
\[
n_3 \geq ha + s_{n_3}x \geq  ha +(y-1)x. 
\]
Similarly, 
\[
s_{n_4} \geq y-1.
\]
The inequality $t_{n_4} \leq h-s_{n_4}$ implies 
\begin{align*}
n_4 & = ha + s_{n_4}x + t_{n_4}y \\
& \leq ha + s_{n_4}x + (h-s_{n_4}) y \\
& = ha + hy - s_{n_4}(y-x) \\
& \leq ha + hy -  (y-1) (y-x).
\end{align*}
Therefore, 
\begin{align*}
n_4 - n_3 
& \leq (ha + hy -  (y-1) (y-x)) - (ha +(y-1)x ) \\
& = hy - y^2 + y. 
\end{align*}
Applying inequality~\eqref{Aristotle:n3n4ineq} gives
\[
N - 2y + 2 \leq n_4 - n_3 \leq hy - y^2 + y 
\]
and so 
\begin{align*}
n_h(A) = N & \leq - y^2 + (h+3)y -2 \\ 
& = \frac{h^2+6h+1}{4} - \left(y-\left(\frac{h+3}{2}\right)\right)^2 \\
& \leq \frac{h^2+6h+1}{4}. 
\end{align*}
Combining the latter with inequality~\eqref{Aristotle:nhA-ineq}, we deduce
\[
n_h(2) = \left[ \frac{h^2+6h+1}{4} \right]
\]
and
\[
n^{\flat}_h(2) \leq n_h(2). 
\]
This completes the proof. 
\end{proof}

\section{Computing \texorpdfstring{$n^{\flat}_h(2)$}{n-flat(h, 2)}  for even \texorpdfstring{$h$}{h}}
In this section we match the upper bound from the previous section for all even $h \neq 4$. We thereby show that $n^{\flat}_h(2) = n_h(2)$ holds for all even $h \neq 4$, and we first consider the exceptional case.

\bt                  \label{Aristotle:theorem:h=4}
For $h=4$, 
\[
n^{\flat}_4(2) = 8 < 10 = n_4(2).
\] 
The extremal bases for $n^{\flat}_4(2)$ are 
\[
\{-3,1,2\}, \quad \{-2,1,2\}, \quad \{-2,1,3\}, \quad \{-1,1,2\}, \quad \{-1,1,4\}.
\] 
\et

\begin{proof} 
From Theorem~\ref{Aristotle:theorem:Stohr} we obtain $n_4(2) = 10$. As for $n^{\flat}_4(2)$, 
let $A = \{a, b, c\}$ be a set of integers with $a < b < c$ and $a < 0$.  
Setting $x=b-a$ and $y = c-a$, we have $A = \{a,a+x,a+y\}$. 
If $[0,5] \subseteq 4A$, then inequalities~\eqref{Aristotle:xya-y} 
and~\eqref{Aristotle:xya-a} from 
Lemma~\ref{Aristotle:lemma:xya} with $h = 4$ give  
\[
1 \leq x < y \leq 5
\]
and
\[
a > -4.
\]
There are $3\binom{5}{2} = 30$ sets $A  = \{a,a+x,a+y\}$ with $a < 0$
that satisfy these conditions.  
One can check either by hand or by computer that $n_4(A) \leq 8$ 
for each of these 30 sets, with equality for the five sets listed.
\end{proof}
 
For the even values of $h \neq 4$, we need to separately consider whether $h \equiv 0 \pmod{4}$ or $h \equiv 2 \pmod{4}$. Let us start with the latter.

\bt 
If $h \equiv 2 \pmod{4}$, then  
\[
n^{\flat}_h(2) = n_h(2).
\]
\et

\begin{proof}
Writing 
\[
h = 4w-2
\]
we have  
\[
n_h(2) = \frac{h^2+6h}{4} = (w+1)h.
\]
Let 
\[
A = \{-w,w,w+1\} = -w + \{0,2w,2w+1\}.
\]
By Lemma~\ref{Aristotle:lemma:3A}, 
\begin{align*}
hA & = -hw + \left\{ 2ws + (2w+1)t: s,t \in \N_0 \text{ and } s+t \leq h \right\} \\
& = -hw + \bigcup_{\ell=0}^h \left\{ 2ws + (2w+1)t: s,t \in \N_0 \text{ and } s+t = \ell  \right\} \\
& = -hw + \bigcup_{\ell=0}^h \left\{ 2w\ell + t : t \in [0,\ell] \right\} \\
& = -hw + \bigcup_{\ell=0}^h [2w\ell , (2w + 1)\ell].
\end{align*}
We have 
\[
2w(\ell  +1) \leq 2w\ell + \ell + 1 
\] 
if and only if 
\[
\ell \geq 2w-1
\]
and so 
\begin{align*}
hA 
&  = \left(-hw + \bigcup_{\ell=0}^{2w-2} [2w\ell, (2w + 1)\ell] \right)
\cup \left(-hw + \bigcup_{\ell=2w-1}^h [2w\ell, (2w + 1)\ell] \right) \\
&  = \left(-hw + \bigcup_{\ell=0}^{2w-2} [2w\ell, (2w + 1)\ell] \right)
\cup \left(-hw + [2w(2w-1), (2w + 1)h] \right) \\
&  = \left(-hw + \bigcup_{\ell=0}^{2w-2} [2w\ell, (2w + 1)\ell] \right)
\cup   [0, (w + 1)h]  \\
& \supseteq [0, n_h(2)],
\end{align*}
and so $n_h(2) = n_h(A) \leq n^{\flat}_h(2) \leq n_h(2)$.  
This completes the proof. 
\end{proof}

For the proof of the remaining even case we introduce the notation
\[
S_{\ell\pmod{m}} = \{n \in S:n\equiv \ell\pmod{m}\}.
\]

\bt               \label{Aristotle:theorem:nh2-upper}
If  $h \geq 8$ and $h  \equiv 0 \pmod{4}$, then 
\[
n^{\flat}_h(2)  = n_h(2).
\]
\et

\begin{proof}
Writing 
\[
h = 4w
\]
we have $w \geq 2$ and 
\[
n_h(2) = \frac{h^2+6h}{4} = w(4w+6).
\]
Let 
\[
A = \{-(w-1),w,w+2\} = -(w-1)+ \{0,2w-1,2w+1\}.
\]
Note that if $h=4$, then $w=1$ and the set $A = \{0,1,3\}$ 
contains no negative integer. 

By Lemma~\ref{Aristotle:lemma:3A}, 
\begin{align*}
hA & = -h(w-1) + \left\{ (2w-1)s + (2w+1)t: s,t \in \N_0 \text{ and } s+t \leq h \right\} \\
& = -h(w-1) + \bigcup_{\ell=0}^h \left\{ (2w-1)s + (2w+1)t: s,t \in \N_0 \text{ and } s+t = \ell  \right\} \\
& =  -h(w-1) + \bigcup_{\ell=0}^h \left\{ (2w-1)\ell + 2t : t \in [0,\ell] \right\} \\
& =  -h(w-1) + \bigcup_{\ell=0}^h   [(2w-1)\ell, (2w+1)\ell]_{\ell\pmod{2}}.
\end{align*}
The set $ [(2w-1)\ell, (2w+1)\ell]_{\ell\pmod{2}}$ is 
a sequence of $\ell+1$ consecutive even integers if $\ell$ is even and 
a sequence of $\ell+1$ consecutive odd integers if $\ell$ is odd. 
We have 
\[
(2w-1)(\ell  +2) \leq (2w+1)\ell + 2 
\] 
if and only if 
\[
2w-2 \leq \ell.
\]
Because $h$ is even, we have 
\begin{align*}
 -h(w-1) +&   \bigcup_{\substack{\ell=2w-2 \\ \ell \equiv 0 \pmod{2}}}^h 
  [(2w-1)\ell, (2w+1)\ell]_{0\pmod{2}} \\
& =  -h(w-1) +  [(2w-1)(2w-2), (2w+1)h]_{\text{even}} \\
& =   [-2(w-1), (w+2)h]_{\text{even}} \\
& \supseteq [-2(w-1), n_h(2)]_{\text{even}}
\end{align*}
and 
\begin{align*}
 -h(w-1) &+  \bigcup_{\substack{\ell=2w-1 \\ \ell \equiv 1 \pmod{2}}}^{h-1} 
  [(2w-1)\ell, (2w+1)\ell]_{1 \pmod{2}}\\
& =  -h(w-1) + [(2w-1)(2w-1), (2w+1)(h-1)]_{\text{odd}} \\
& = [1, -h(w-1)+(2w+1)(h-1)]_{\text{odd}} \\
& =    [1,n_h(2)-1]_{\text{odd}}. 
\end{align*}
Therefore,
\[
hA \supseteq [0,n_h(2)]
\]
and $n_h(2) \geq n^{\flat}_h(2) \geq n_h(A)$. 
This completes the proof. 
\end{proof}

\section{Bounding \texorpdfstring{$n^{\flat}_h(2)$}{n-flat(h, 2)} for odd \texorpdfstring{$h$}{h}}

\bt 
For all positive odd integers $h$,  
\[
n^{\flat}_h(2) < n_h(2).
\]
\et

\begin{proof} 
It is straightforward to check that  
\[
n^{\flat}_1(2) = 1 < 2 = n_1(2) 
\]
and that $\{0,1,2\}$ is  the unique extremal set for $n_1(2)$.   
The extremal sets for $n^{\flat}_1(2)$ are the sets $\{a, 0,1\}$ for 
all  integers $a < 0$.  

Let $h = 2m+1 \geq 3$.  
By Theorem~\ref{Aristotle:theorem:new}, 
\[
n^{\flat}_h(2)  \leq n_h(2) =  \frac{h^2 + 6h +1}{4} = m^2 + 4m + 2 = N
\]
and $h+1 \leq N$.  
If $n^{\flat}_h(2) = n_h(2)$, then there is a set $A = \{a,b,c\}$ with 
$a < b < c$ and $a < 0$ such that 
\[
[0,h+1] \subseteq [0, N] \subseteq hA. 
\]
With $x = b-a$ and $y = c-a$, we define $n_3$ and $n_4$ in the exact same way as in the proof of Theorem~\ref{Aristotle:theorem:new}.
We then still have
\[
y-1 \ge n_3 \ge ha + (y-1)x,
\]
which, on the one hand, gives
\[
ha \le (y-1)(1-x) = U.
\]
On the other hand, we also still have
\[
N - y + 1 \le n_4 \le ha + hy - (y-1)(y-x),
\]
implying
\[
ha \ge N - hy + (y-1)(y - x - 1) = L.
\]
In other words,
\[
0 \geq L-U = N - hy + (y-1)(y-2).
\]
Recalling that $h = 2m+1$ and $N = m^2 + 4m + 2$, we obtain 
\begin{align*}
0 & \geq L - U = N - hy + (y-1)(y-2)  \\
& = (m^2+4m+2) -(2m+1)y+y^2-3y+2 \\
& = (m-y)^2 +4(m-y)+4\\
& = (m-y+2)^2 \\ 
& \geq 0. 
\end{align*}  
\begin{samepage}
Therefore,
\[
y = m + 2 
\]
and 
\[
ha = L = U,
\]
so that
\[
ha = (y-1)(1-x) = (m+1)(1-x) < 0. 
\]
\end{samepage}
It follows that $1-x < 0$ and so 
\[
1 \leq x -1 \leq y - 2 = m.
\]
Also, 
\[
2ha = 2(m+1)(1-x) = (h+1)(1-x) = h(1-x)+1-x
\]
implies 
\[
(2m+1)(2a+x-1) = h(2a+x-1) = -(x-1)
\]
and so $2m+1$ divides the integer $x-1$ in $[1,m]$, which is absurd.
Therefore, $n^{\flat}_h(2) < n_h(2)$, completing the proof. 
\end{proof}

\section{AI usage}
The first proofs of the comparison results for the extremal functions $n_h(2)$
and $n^{\flat}_h(2)$ were autonomously discovered and formally verified in Lean
by \emph{Aristotle}, a formal reasoning agent developed by
\href{https://www.harmonic.fun/}{Harmonic}~\cite{Aristotle}.
WvD and PM used Aristotle to obtain these original proofs.
Aristotle is available at
\href{https://aristotle.harmonic.fun/}{\nolinkurl{aristotle.harmonic.fun}}.

MBN studied Aristotle's arguments and used its conversational interface to
clarify the reasoning. He then developed the new proofs presented in this
paper, incorporating some of Aristotle's methods, and wrote the original
manuscript.

\end{document}